%&amstex
\magnification=1200
\NoBlackBoxes
\documentstyle{amsppt}
%\TagsOnRight
\topmatter
\title  Hankel Multipliers And Transplantation Operators
\endtitle
\rightheadtext{Hankel Multipliers And Transplantation}
\author Krzysztof Stempak and Walter Trebels
\endauthor
\affil Instytut Matematyczny Polskiej Akademii Nauk\\
and\\ Technische Hochschule Darmstadt
\endaffil
\address
Instytut  Matematyczny,
Polska Akademia Nauk,
ul. Kopernika 18,
51-617 Wroc\l aw,
Poland.
\endaddress
\address
Fachbereich Mathematik,
TH Darmstadt,
Schlossgartenstr. 7,
D-64289 Darmstadt,
Germany.
\endaddress
\dedicatory Dedicated to Professor Satoru Igari on the occasion
of
his 60th birthday 
-- 10-1-96 version
\enddedicatory
\thanks Research of the first author was supported in part by KBN 
grant \# 2 PO3A 030 09
\endthanks
\abstract
Connections between Hankel transforms of different order for 
$L^p$-functions are examined. Well known are the results of Guy
[Guy] and 
Schindler [Sch]. Further relations result from projection
formulae  for 
Bessel functions of different order. Consequences for Hankel
multipliers are 
exhibited and implications for radial Fourier multipliers on
Euclidean spaces 
of different dimensions indicated.
\endabstract
\endtopmatter
%\eject
\document
\subhead
1. Introduction
\endsubhead
It is well known that harmonic analysis of radial functions on
the Euclidean 
space $\bold{R}^{n}$, $n\ge1$,
reduces to studying appropriate function spaces on the half-line
equipped with the measure $x^{n-1}dx$. The Fourier transform is
then replaced
by the modified Hankel transform of order $\frac{n-2}{2}$. The
aim of this 
paper is to show, among others, that also studying the
non-modified 
Hankel transform of an arbitrary order $\nu\ge-1/2$ within an
appropriate 
weighted
setting leads to corresponding results for Fourier transform on
radial 
functions. This is seen, for instance, in Section 2 where we
discuss multiplier
results for the modified Hankel transform. It occurs that they
are closely
related to two transference theorems of Rubio de Francia for
Fourier transform
on Euclidean spaces.

Given $\nu\ge-1/2$ and $f$, an integrable function on $\bold
{R}_{+}=(0,\infty)$, its (non-modified) Hankel transform is
defined by
$$
\Cal H_{\nu}f(x)=\int_{0}^{\infty}(xy)^{1/2}J_{\nu}(xy)f(y)dy,
\qquad
x>0. \tag1.1
$$ 
Here $J_{\nu}(x)$ denotes the Bessel function of the first kind
of
order $\nu$, [Sz, (1.17.1)]. For $\nu=-1/2$ or $\nu=1/2$ one
recovers
the cosine and sine transforms on the half-line.
The modified Hankel transform is given by
$$
H_{\nu}f(x)=\int_0^\infty \frac{J_{\nu}(xy)}{(xy)^{\nu}} f(y)
y^{2\nu+1}dy, \qquad x>0. \tag1.2
$$
Due to the estimates on Bessel function
$$
J_{\nu}(x)=O(x^{\nu}),\qquad J_{\nu}(x)=O(x^{-1/2}), \tag1.3
$$
valid for $x\to0$ and $x\rightarrow\infty$ correspondingly,
$H_{\nu}f$ is
well-defined for every function $f$ in
$L^{p}(\bold{R}_{+},x^{2\nu+1}dx)$,
$1\le p<\frac{4(\nu+1)}{2\nu+3}$.
Clearly both transforms are related to each other by
$$
\Cal
H_{\nu}f(x)=x^{\nu+1/2}H_{\nu}((\cdot)^{-(\nu+1/2)}f(\cdot))(x)\tag1.4
$$
whenever $f$ is an integrable function on $\bold{R}_+$ and, for
instance, $\int_0^\infty|f(x)|x^{\nu+1/2}dx<\infty$.
Moreover, for $\nu\ge-1/2$ the inversion formulae
$$
f(y)=\int_{0}^{\infty}(xy)^{1/2}J_\nu(xy)\Cal H_{\nu}f(x)dx
\tag1.5
$$
and
$$
f(y)=\int_{0}^{\infty}\frac{J_{\nu}(xy)}{(xy)^{\nu}}H_{\nu}f(x)
x^{2\nu+1}dx \tag1.6
$$
hold: (1.5) holds, for instance, for every $C^1$ function $f\in
L^1(\bold
{R}_+, dx)$ with $\Cal H_\nu f\in L^1(\bold{R}_+, dx)$; (1.6)
holds for every 
$C^1$ function $f\in L^1(\bold{R}_+, x^{2\nu+1}dx)$ with 
$H_{\nu}f\in L^1(\bold{R}_+, x^{2\nu+1}dx)$, cf\. [W, p\.456].
More can be said: $H_\nu$ is a bijection on  
$\Cal S(\bold{R}_{+})$, the space of infinitely
differentiable even functions on $\bold{R}$ with rapidly
decreasing
derivatives, while $\Cal H_\nu$ is a bijection on the Zemanian
space
$\Cal Z_\nu$ of all $C^\infty$ functions $f$ on $\bold{R}_+$ for
which
the quantity
$$
\sup_{x>0}|x^n\bigg(\frac1x\frac d{dx}\bigg)^k(x^{-\nu-1/2}f(x))|
$$
is finite for every $n, k\in\bold{N}=\{0,1,2,\dots\}$, see [Z1,
Z2].
Note at this point that
$C^{\infty}_o=C^{\infty}_o(\bold{R}_{+})$,
the space of compactly supported
$C^\infty$ functions on $\bold{R}_+$, is contained in every $\Cal
Z_\nu$.
The kernels $\varphi^{\nu}_{x}(y)=(xy)^{1/2}J_{\nu}(xy)$, $y>0$,
appearing in (1.1) satisfy
$$
\bigg(\frac{d^{2}}{dy^{2}}+\frac{1/4-\nu^{2}}{y^{2}}\bigg)
\varphi^{\nu}_{x}(y)=-x^{2}
\varphi^{\nu}_{x}(y), \qquad x>0.\tag1.7
$$
while the kernels 
$\phi^{\nu}_{x}(y)=(xy)^{-\nu}J_{\nu}(xy)$, $y>0$,
appearing in (1.2) fulfil
$$
\bigg(\frac{d^{2}}{dy^{2}}+\frac{2\nu+1}{y}\frac{d}{dy}\bigg)
\phi^{\nu}_{x}(y)=
-x^{2}\phi^{\nu}_{x}(y),\qquad x>0. \tag1.8
$$
The differential operators on the left sides of (1.7) and (1.8)
are
symmetric in $L^{2}(\bold{R}_{+}, dx)$ and $L^{2}(\bold
{R}_{+},x^{2\nu+1}dx)$ correspondingly.

\noindent As usual we use $C$ or $c$ with or without subscripts
as a 
constant which is not necessarily the same at each occurence.
\subhead
2. Hankel multipliers
\endsubhead
In this section we fix $\nu\ge-1/2$ and consider weighted
Lebesgue
spaces on $\bold{R}_{+}$ with respect to the Lebesgue measure
$dx$ on
one occasion and the measure
$$
dm_{\nu}(x)=x^{2\nu+1}dx
$$
on another one. Hence, in what follows we use the notation
$$
||f||_{p,\alpha}=\bigg(\int_{0}^{\infty}|f(x)|^{p}x^{\alpha}dx\bigg)^{1/p}
$$
and
$$
||f||_{L^{p}(x^{\alpha}dm_{\nu})}=\bigg(\int_{0}^{\infty}|f(x)|^{p}x^{\alpha}d
m_{\nu}(x)\bigg)^{1/p}
$$
for $1\le p<\infty$ with usual modification when $p=\infty$. By
$L^{p,\alpha}(dx)$ and $L^{p,\alpha}(dm_{\nu})$ we denote the
weighted Lebesgue spaces of functions for which the above
quantities
are finite. If $\alpha=0$ we write $L^p$ instead of $L^{p,0}$.
By $\Cal M^{\nu,\alpha}_{p}$, $M^{\nu,\alpha}_{p}$ we
denote the spaces of weighted $p$-multipliers for the Hankel and
modified Hankel transform.
Thus, a bounded measurable function $m(x)$
on $\bold{R}_{+}$ is in $\Cal M^{\nu,\alpha}_{p}$ provided
$$
||\Cal H_{\nu}(m\cdot \Cal H_{\nu}f)||_{p,\alpha}\le
C||f||_{p,\alpha}\,,
$$
where C is a constant independent of $f$ in $\Cal
H_\nu(C^\infty_o)$,
the image of $C^\infty_o$ under the action of $\Cal H_\nu$. 
The least constant C for
which the above inequality holds is called the multiplier norm of
$m$.
Similar definition is for the multiplier space
$M^{\nu,\alpha}_{p}$, 
now with the norm $||\cdot||_{L^p(x^\alpha dm_\nu)}$ in use, and
here
$H_\nu(C^\infty_o)$ is the testing function space.

We postpone to Section~4
the proof of the fact that $H_\nu(C^\infty_o)$ is dense in 
$L^{p,\alpha}(dx)$ if $1<p<\infty$ and $\alpha>-1$ while 
$\Cal H_\nu(C^\infty_o)$ is dense in $L^{p,\alpha}(dx)$ if
$1<p<\infty$
and $\alpha>-1-p(\nu+1/2)$ (the case $p=1$ requires additional
assumptions).
This is the contents of Theorem~4.7 and Corollary~4.8.

The following is Guy's transplantation theorem for the Hankel
transform (cf. also [Sch] for an alternative proof).

\proclaim{Theorem ([Guy, Lemma 8C])}
Let $\mu,\nu\ge-1/2$, $1<p<\infty$ and $-1<\alpha<p-1$. Then
$$
C^{-1}||\Cal H_{\nu}f||_{p,\alpha}\le||\Cal H_{\mu}f||_{p,\alpha}
\le C||\Cal H_{\nu}f||_{p,\alpha}
$$
with $C=C(\mu,\nu,p,\alpha)$ independent of $f\in L^1(\bold{R}_+,
dx)$.
\endproclaim
As an immediate consequence one obtains
\proclaim{Corollary~2.1} 
Let $\mu,\nu\ge-1/2$, $1<p<\infty$ and $-1<\alpha<p-1$. Then
$$
\Cal M^{\nu,\alpha}_{p}=\Cal M^{\mu,\alpha}_{p}. \tag2.1
$$ 
\endproclaim

\demo{Proof}
Assuming $m$ is in  $\Cal M^{\nu,\alpha}_{p}$ and 
$f$ is in $\Cal H_\mu(C^\infty_o)$ and using
the fact that $\Cal H_{\mu}f\in C^\infty_o\subset L^1(dx)$, hence

$\Cal H_\nu\Cal H_\mu f\in \Cal H_\nu(C^\infty_o)$, we write
$$
\align
||\Cal H_{\mu}(m\cdot \Cal H_{\mu}f)||_{p,\alpha}&\le
C||\Cal H_{\nu}(m\cdot\Cal H_{\nu}(\Cal H_{\nu}\Cal
H_{\mu}f))||_{p,\alpha}\\
&\le CC_{\nu,m}||\Cal H_{\nu}\Cal H_{\mu}f||_{p,\alpha}\\
&\le C^{2}C_{\nu,m}||f||_{p,\alpha},
\endalign
$$
where $C_{\nu,m}$ denotes the operator
norm of the multiplier $m\in \Cal M^{\nu,\alpha}_{p}$.  Thus
$\Cal M^{\nu,\alpha}_{p} \subset\Cal M^{\mu,\alpha}_{p}$. 
Analogously
the opposite inclusion follows.
\enddemo
\proclaim{Corollary~2.2}
Let $\mu,\nu\ge-1/2$, $1<p<\infty$. Assume further
that $-1<\beta+(\nu+1/2)(2-p)<p-1$ and denote $\beta^{*}=\beta+
(\nu-\mu)(2-p)$. Then
$$
M^{\nu,\beta}_{p}=M^{\mu,\beta^{*}}_{p}. \tag2.2
$$
\endproclaim
\demo{Proof}
The identity (1.4), the fact that
$x^{\nu+1/2}C^\infty_o=C^\infty_o$ 
and the definition of multiplier spaces immediately give
$$
M^{\nu,\beta}_{p}=\Cal
M^{\nu,\beta+p(2\nu+1)(\frac1p-\frac12)}_{p}
$$
and then (2.1) produces (2.2).
\enddemo
In particular, (2.2) for $p=2$ gives
\proclaim{Corollary~2.3}
Let $\mu,\nu\ge-1/2$ and $-1<\beta<1$. Then
$$
M^{\nu,\beta}_{2}=M^{\mu,\beta}_{2}. \tag2.3
$$
\endproclaim
In some sense (2.3) may be viewed as a ``radial'' generalization
of
Rubio de Francia transference result, [RdF, Theorem~2.1], which
claims that
given $-1<\beta<1$ and $m\in L^\infty(\bold{R}_+)$, being a
Fourier multiplier
by $m(|x|)$ on $L^2(\bold{R}\,, |x|^\beta dx)$, implies
$m(||x||)$ to be a
Fourier multiplier on $L^2(\bold{R}^n, ||x||^\beta dx)$, $n\ge 2$
($||\cdot||$ denotes here the Euclidean norm in appropriate
Euclidean
space,  $dx$  the Lebesgue measure). When restricted to the space
of
radial functions on which the multiplier acts, (2.3) claims that
the
opposite implication also holds. In particular this implies that
one
can jump between Euclidean
spaces of arbitrary dimension in contrast to the preceding
situation where
a jump was allowed only between $\bold{R}$ and $\bold{R}^n$.

Let $T_{R}$, $R>0$, denote the multiplier operator corresponding
to
the characteristic function  of the interval $(0,R)$. By a
homogeneity argument, for every $\mu, p$ and $\alpha$, the
operator
norms of $T_{R}$ as members of $\Cal M^{\mu,\alpha}_{p}$ or
$M^{\mu,\alpha}_{p}$ are independent of $R>0$. Hence, in what
follows
consider $T=T_{1}$ only. Hirschman's ``one-dimensional'' weighted
multiplier result, [Hi1], says that $T\in\Cal M^{-1/2,\alpha}_{p}
=M^{-1/2,\alpha}_{p}$ for every $1<p<\infty$ and $-1<\alpha<p-1$.
Thus (2.1) further gives that $T\in\Cal M^{\nu,\alpha}_{p}$ for
every
$\nu>-1/2$, $1<p<\infty$ and $-1<\alpha<p-1$ which, for
$\alpha=0$
was proved by Wing, [Wi]. Similarly, Herz' result, [He], which
says
that $T\in M^{\nu,0}_p$ for $\nu>-1/2$ provided
$$
\frac{4(\nu+1)}{2\nu+3}<p<\frac{4(\nu+1)}{2\nu+1}
$$
may be recovered from (2.2) and just mentioned Hirschman's result

The next corollary may be considered as a ``radial'' extension of
another result due to Hirschman [Hi2].
\proclaim{Corollary~2.4}
Let $\mu\ge-1/2$ and $-1<\beta<1$. Then $T\in M^{\mu,\beta}_{2}$.
\endproclaim
\demo{Proof}
Combine (2.3) and the fact that $T\in M^{-1/2,\beta}_{2}$.
\enddemo
As already mentioned, Schindler  gave an alternative proof of
Guy's result. Besides, in the special case $\mu=\nu+2k$,
$k=1,2,\dots$, her approach allowed to take into account both
endpoints $p=1$ and $p=\infty$ and, in addition, to obtain a
range of
$\alpha$'s different than the $A_p$ range from Guy's theorem.
\proclaim{Theorem ([Sch, Theorems 3 and 4])}
Let $\nu\ge-1/2$, $1\le p<\infty$, $k=1,2,\dots$ and 
$-p(\nu+1/2)<\alpha<p(\nu+1/2)$. Then
$$
C^{-1}||\Cal H_\nu f||_{p,\alpha}\le||\Cal
H_{\nu+2k}f||_{p,\alpha}
\le C||\Cal H_{\nu}f||_{p,\alpha}
$$
with $C=C(\nu,k,p,\alpha)$ independent of $f\in L^1(\bold{R}_+,
dx)$.
\endproclaim
In consequence, analogous to (2.1) is now the identity
$$
\Cal M^{\nu,\alpha}_{p}=\Cal M^{\nu+2k,\alpha}_{p},\tag2.4
$$ 
where $\nu\ge-1/2$, $1\le p<\infty$, $k=1,2,\dots$ and 
$-p(\nu+1/2)<\alpha<p(\nu+1/2)$.
Moreover, in the case $p>1$, comparing the above hypotheses with
those from Corollary~2.1
allows further to enlarge the range of $\alpha$'s for which (2.4)
holds to
$$
-\max\{p(\nu+1/2),1\}<\alpha<\max\{p(\nu+1/2),p-1\}.
$$
Similarly, analogous to (2.2) is
$$
M^{\nu,\beta}_{p}=M^{\nu+2k,\beta^*}_{p},\tag2.5
$$ 
where $\nu\ge-1/2$, $1<p<\infty$, $k=1,2,\dots$, 
$\beta^*=\beta-2k(2-p)$  and 
$$
-\max\{p(\nu+1/2),1\}<\beta+p(2\nu+1)(\frac1p-\frac{1}{2})<\max\{
p(\nu+1/2),
p-1\}.
$$
In particular, (2.5) for $p=2$ gives
\proclaim{Corollary~2.5}
Let $\nu\ge-1/2$,  $k=1,2,\dots$ and 
$-\max\{2\nu+1,1\}<\beta<\max\{2\nu+1,1\}$.
Then 
$$
M^{\nu,\beta}_{2}=M^{\nu+2k,\beta}_{2}.
$$ 
\endproclaim
The above  stands in a relationship with another Rubio de Francia
transference result, [RdF, Theorem 2.2], in the same
way as (2.3) ``generalizes" [RdF, Theorem 2.1]. This result says
that given
$m\in L^\infty(\bold{R}_+)$ and $w(s)$, a nonnegative measurable
function
on $\bold{R}_+$, if $m(||x||)$ is a Fourier multiplier on
$L^2(\bold{R}^n,
w(||x||)dx)$ for some $n\ge2$ then $m(||x||)$ is a Fourier
multiplier
 on $L^2(\bold{R}^{n+2k}, w(||x||)dx)$ for any $k=1,2,\dots$.
Corollary~2.5 says that when restricted to indicated power
weights and
spaces of radial functions on which multipliers act, the converse
also
holds provided the difference in Euclidean dimensions is a
multiple
of 4 (enlarging $\nu$ by 2 changes the Euclidean dimension by 4).
Speaking less precisely this means that under appropriately
modified assumptions we can exchange radial Fourier multipliers,
in both
directions, between Euclidean spaces whose difference in
dimensions is
a multiple of 4 (one-dimensional situation is now included!).

\subhead
3. Weighted estimates for the transference operators
\endsubhead
Throughout this section all the functions we are dealing with are
assumed to be in $C^{\infty}_o$.
Let 
$$
L_\nu=-\bigg(\frac{d^{2}}{dy^{2}}+\frac{2\nu+1}{y}\frac{d}{dy}\bigg),
 \qquad\nu\ge-1/2\,,
$$
be the differential operator appearing in (1.8). Clearly
$$
H_\nu(L_\nu f)(y)=y^2H_\nu f(y).
$$
Hence, in terms of the modified Hankel transform,
$$
H_\nu(L_\nu^\delta f)(y)=y^{2\delta}H_\nu f(y) \tag 3.1
$$
is a well motivated definition of $L_\nu^\delta$, the
$\delta$-fractional power of $L_\nu$. Rewriting in terms of the
modified Hankel transform the inequality
$$
||\Cal H_\mu\Cal H_\nu f||_{p,\alpha}\le C||f||_{p,\alpha}\,\,,
$$
which follows from Guy's transplantation theorem, gives
$$
\bigg(\int_0^\infty|H_\mu H_\nu h(x)|^px^\gamma dx\bigg)^{1/p}
\le C\bigg(\int_0^\infty|L_\nu^{\frac{\mu-\nu}2}h(x)|^px^\delta
dx\bigg)^{1/p},
$$
where $\gamma=p(\mu+1/2)+\alpha$, $\delta=p(\nu+1/2)+\alpha$ for
$1<p<\infty$, $-1<\alpha<p-1$ and $\mu, \nu\ge-1/2$.

In this section we prove two weighted $L^p$--$L^q$ inequalities
for the
transplantation operator
$$
T_\nu^\mu=H_\mu\circ H_\nu .
$$
This is achieved first by using appropriately chosen integral
formulae for
Bessel functions that generate nice representations of
$T^{\mu}_{\nu}$
%$$
%T_\nu^\mu
%g(x)=c_{\nu,\mu}\cdot\frac1{x^{2\mu}}\int_0^x(x^2-y^2)^{\mu-\nu-1}
%L_\nu^{\mu-\nu} g(y)y^{2\nu+1}dy
%$$
%and
%$$
%T_\nu^\mu g(x)=c_{\nu,\mu}\int_x^\infty(y^2-x^2)^{\nu-\mu-1}
%yg(y)dy
%$$
(with necessary restrictions on $\nu$ and $\mu$) and then
applying
 some weighted norm inequalities for the Riemann-Liouville
and Weyl fractional integral operators.
\proclaim{Theorem 3.1}
Let $-1/2\le\nu<\mu$, $1<p<q<\infty$, $p(\mu+1)\ge1$ and 
$$
\frac{\mu+1}q=\frac{\nu+1}p\,. \tag 3.2
$$
Then
$$
\bigg(\int_0^\infty|T_\nu^\mu g(x)|^qx^{2\mu+1}dx\bigg)^{1/q}\le
C\bigg(\int_0^\infty|L_\nu^{\mu-\nu}g(x)|^px^{2\nu+1}dx\bigg)^{1/
p}.
\tag3.3 
$$
\endproclaim            

\demo{Proof} 
By using a homogeneity argument it is easy to see that (3.2) is
necessary
for (3.3) to hold. If $\rho_rg(x)=\frac 1r g(\frac xr)$, $r>0$,
then
$H_\nu(\rho_rg)(x)=r^{2\nu+1}H_\nu g(rx)$
and
$$
L^{\mu-\nu}_\nu(\rho_rg)(x)=r^{2(\nu-\mu)}\rho_r(L^{\mu-\nu}_\nu
g)(x)\,.
$$
Hence, considering (3.3) with $\rho_rg$ in place of $g$ and
allowing $r$
to be small and large gives (3.2).

Evaluating the formula [EMOT, 8.5\,(33)]  at $y=1$ and writing
$\mu$
in place of $\mu+\nu+1$ produces
$$
\frac{J_\mu(a)}{a^\mu}=c_{\nu,\mu}\frac1{a^{2\mu}}
\int_0^a (a^2-t^2)^{\mu-\nu-1}\frac{J_\nu(t)}{t^\nu}t^{2\nu+1}dt
$$
for $-1<\nu<\mu$ and $a>0$.
Hence, a change of variable and Fubini's theorem give
$$
\align
T_\nu^\mu g(x)
&=\int_0^\infty H_\nu
g(y)\frac{J_{\mu}(xy)}{(xy)^{\mu}}y^{2\mu+1}dy\\
&=c_{\nu,\mu}\int_0^\infty H_\nu g(y) \frac1{(xy)^{2\mu}}
\int_0^{xy}
((xy)^2-t^2)^{\mu-\nu-1}\frac{J_\nu(t)}{t^\nu}t^{2\nu+1}dt
y^{2\mu+1}dy\\
&=c_{\nu,\mu}\frac1{x^{2\mu}}\int_0^x(x^2-u^2)^{\mu-\nu-1}\int_0^
\infty
y^{2(\mu-\nu)}H_\nu
g(y)\frac{J_\nu(uy)}{(uy)^\nu}\,y^{2\nu+1}dyu^{2\nu+1}du\,.
\endalign
$$
Note that an application of Fubini's theorem is possible due to
the 
boundedness of $J_\nu(s)/s^\nu$ on $(0,\infty)$, integrability of

$y^{2\mu}H_\nu g(y)$ on $(0,\infty)$ (with respect to the
Lebesgue measure) and the 
assumption $\mu-\nu>0$.

Taking into account (3.1) and the inversion formula for the
modified
Hankel transform (recall that $\nu\ge-1/2$) we arrive at
$$
T_\nu^\mu
g(x)=c_{\nu,\mu}\cdot\frac1{x^{2\mu}}\int_0^x(x^2-y^2)^{\mu-\nu-1
}
L_\nu^{\mu-\nu} g(y)y^{2\nu+1}dy\,.
$$
What we now need is the inequality
$$
\align
\bigg(\int_0^\infty\bigg(\frac1{x^{2\mu}}\int_0^x(x^2-y^2&)^{\mu-
\nu-1}
G(y)y^{2\nu+1}dy\bigg)^qx^{2\mu+1}dx\bigg)^{1/q}\\
&\le C
\bigg(\int_0^\infty G(x)^px^{2\nu+1}dx\bigg)^{1/p},
\endalign
$$
say, for all nonnegative functions $G$. Elementary variable
changes show 
that the above inequality is equivalent to
$$
\align
\bigg(\int_0^\infty\bigg(\int_0^t(t-&s)^{\mu-\nu-1}h(s)ds\bigg)^q
t^{\mu(1-q)}dt\bigg)^{1/q}\tag3.4\\
&\le C\bigg(\int_0^\infty h(t)^p
t^{\nu(1-p)}dt\bigg)^{1/p},
\endalign
$$
$h-$ nonnegative. We use the following criterion for $L^p-L^q$
weighted
estimates for the Riemann-Liouville 
$$
I^{\alpha}_+h(t)=\frac1{\Gamma(\alpha)}\int_0^t(t-s)^{\alpha-1}h(
s)ds
$$
and Weyl
$$
I^{\alpha}_-h(t)=\frac1{\Gamma(\alpha)}\int_t^\infty
(s-t)^{\alpha-1}h(s)ds
$$
fractional integral operators.
\proclaim{Lemma ([SKM, Theorem 5.4])}
Let $\alpha>0$, $p\ge1$ and $p\le q\le p/(1-p\alpha)$ when $1\le
p<
1/\alpha$ (if $p=1$ then the right endpoint $p/(1-p\alpha)$ is
excluded)
or $p\le q<\infty$ when $p\ge1/\alpha$. Suppose also that
$-\infty<N<
\infty$ and $M<p-1$ when we consider $I^{\alpha}_+$  or $M>\alpha
p-1$ when
$I^\alpha_-$ is taken into account and
$$
\frac{N+1}q=\frac{M+1}p-\alpha. \tag 3.5
$$
Then
$$
\bigg(\int_0^\infty|I_\pm^\alpha h(t)|^qt^Ndt\bigg)^{1/q}
\le C\bigg(\int_0^\infty|h(t)|^pt^Mdt\bigg)^{1/p}.\tag3.6
$$
\endproclaim
To see that the lemma  gives (3.4) and thus (3.3) assume that the
hypotheses 
of Theorem 3.1 are satisfied and take $\alpha=\mu-\nu$,
$N=\mu(1-q)$,
$M=\nu(1-p)$. Clearly (3.2) gives (3.5) and $M<p-1$ holds
provided
$\nu>-1$. Moreover, if $1<p<1/\alpha$ then the condition
$p(\mu+1)\ge1$
implies $q\le p/(1-p\alpha)$. Thus, the conclusion of the above
lemma, (3.6),
holds for the operator $I^\alpha_+$ and, in consequence, (3.4) is
valid. 
This concludes the proof of Theorem 3.1.
\enddemo
\proclaim{Theorem 3.2}
Let $-1/2\le\mu<\nu$, $1\le p\le q<\infty$ and $p\mu+1\ge0$.
Suppose also that
$$
\nu+\frac{\mu+1}q=\mu+\frac{\nu+1}p\,.\tag 3.7
$$
Then
$$
\bigg(\int_0^\infty|T_\nu^\mu g(x)|^qx^{2\mu+1}dx\bigg)^{1/q}\le
C\bigg(\int_0^\infty|g(x)|^px^{2\nu+1}dx\bigg)^{1/p}.
\tag3.8
$$
\endproclaim
\demo{Proof}
A homogeneity argument similar to that from the beginning of the
proof
of Theorem 3.1 also shows that (3.7) is necessary for (3.8) to
hold.

Evaluating the formula [EMOT, 8.5\,(32)]  at $y=1$ and writing
$\nu-\mu-1$
in place of $\mu$ produces
$$
\frac{J_\mu(a)}{a^\mu}=c_{\nu,\mu}
\int_a^\infty t^{1-\nu}(t^2-a^2)^{\nu-\mu-1}J_\nu(t)dt
$$
for arbitrary $\nu, \mu$ satisfying $\Re\mu<\Re\nu<2\Re\mu+3/2$
and $a>0$.
Note at this point that only with the stronger assumption 
$\Re\mu<\Re\nu<2\Re\mu+1/2$, $\Re\mu>-1/2$, is the above integral
Lebesgue 
integrable; 
otherwise it converges in the Riemann sense.
Hence, considering first the case $\mu>-1/2$, for real $\nu,\mu$
that 
satisfy $-1/2<\mu<\nu<2\mu+1/2$,
a change of variable and Fubini's theorem give
$$
\align
T_\nu^\mu g(x)
&=\int_0^\infty H_\nu
g(y)\frac{J_{\mu}(xy)}{(xy)^{\mu}}y^{2\mu+1}dy\\
&=c_{\nu,\mu}\int_0^\infty H_\nu g(y)
\int_{xy}^\infty t^{1-\nu}(t^2-(xy)^2)^{\nu-\mu-1}J_\nu(t)dt
\cdot y^{2\mu+1}dy\\
&=c_{\nu,\mu}\int_x^\infty
u^{1-\nu}(u^2-x^2)^{\nu-\mu-1}\int_0^\infty
H_\nu g(y)J_\nu(uy)y^{2\nu+1}dy\,du\\
&=c_{\nu,\mu}\int_x^\infty u(u^2-x^2)^{\nu-\mu-1}\int_0^\infty
H_\nu g(y)\frac{J_\nu(uy)}{(uy)^\nu}y^{2\nu+1}dy\,du
\,.
\endalign
$$
An application of Fubini's theorem is allowed at this point since
$|J_\nu(s)|\le Cs^{-1/2}$ on $(0,\infty)$, the function
$y^{\nu+1/2}
H_\nu g(y)$ is integrable on $(0,\infty)$ and
$u^{1/2-\nu}(u^2-x^2)^{\nu-\mu-1}$ is integrable on $(x,\infty)$
(both
with respect to the Lebesgue measure).

Since $\nu>-1/2$ the inversion formula for the modified Hankel
transform gives
$$
T_\nu^\mu g(x)=c_{\nu,\mu}\int_x^\infty(y^2-x^2)^{\nu-\mu-1}
yg(y)dy\,.\tag3.9
$$

It is now easy to see that the above argument, in particular the
inversion
formula, remains valid for complex $\nu,\mu$ satisfying
$-1/2<\Re\mu<\Re\nu<2\Re\mu+1/2$. Moreover, for any fixed 
$g\in\Cal S(\bold{R}_+)$, $x\in\bold{R}_+$ and $\mu$ with
$\Re\mu>-1/2$,
both sides of (3.9) are analytic functions of the complex
variable $\nu$,
$\Re\nu>\Re\mu$ (analyticity of the coefficient $c_{\nu,\mu}$
follows
from its explicit form, cf\. [EMOT, 8.5(32)]). Hence, the
validity of (3.9)
is implied for every real $\nu,\mu$ with  $-1/2<\mu<\nu$, by an
analytic
continuation method. For the remaining case $\mu=-1/2$ we need
some parameter
interval to
do analytic continuation. Fortunately, the interchange of
integration we did
above is also allowed in the range $-1/2<\nu<1/2$ when
$\mu=-1/2$: for this
more subtle argument see the remarks in [GT] following the proof
of 
[GT, (1.5)]. 

To prove (3.8) we now need the inequality
$$
\bigg(\int_0^\infty\bigg(\int_x^\infty(y^2-x^2)^{\nu-\mu-1}yG(y)d
y
\bigg)^qx^{2\mu+1}dx\bigg)^{1/q}\le C\bigg(\int_0^\infty
G(x)^px^{2\nu+1}
dx\bigg)^{1/p}, 
$$
which, after an elementary change of variable, turns out to be
equivalent
to 
$$
\align
\bigg(\int_0^\infty\bigg(\int_t^\infty(s-t)^{\nu-\mu-1}&h(s)ds\bigg)^qt^{\mu}
dt\bigg)^{1/q} \tag3.10\\
&\le C\bigg(\int_0^\infty h(s)^ps^{\nu}ds\bigg)^{1/p}.
\endalign
$$
Both, $G$ and $h$ are assumed to be nonnegative functions. (3.10)
is now a
consequence of (3.6) for the operator $I^{\alpha}_-$ if we take 
$\alpha=\nu-\mu$, $N=\mu$ and $M=\nu$. Indeed, assuming the
hypotheses
of Theorem 3.2 to be satisfied it is easy to see that (3.7)
forces the
condition $p<(\nu+1)/(\nu-\mu)$ which is nothing else but
$M>\alpha p-1$.
Further, in the case $p<1/(\nu-\mu)$ the assumption $p\mu+1\ge0$
guarantees
$q\le p/(1-p\alpha)$ to hold.
This finishes the proof of Theorem 3.2.
\enddemo

\noindent

Theorems~3\.1 and 3\.2 also have  applications to radial Fourier
multipliers. 
Setting $H_\nu g(y) =m(y)$ Theorem~3\.1 claims (under the 
relevant conditions on the parameters) that 
$$
||H_\mu m||_{L^q(dm_\mu)} \le C||H_\nu (y^{2(\mu -\nu)}m)|| 
_{L^p(dm_\nu)} \, \tag3.11
$$
provided that $H_\nu m\in C^\infty_o(0,\infty)$.
Denote now by $L^p_{rad}({\bold R}^n)$ the set of radial
$L^p$-functions on ${\bold R}^n$, $f(x)=f_o(||x||)$, with
standard 
$L^p({\bold R}^n,dx)$-norm 
and by $[L^p_{rad}({\bold R}^n)]\,\,\widehat{}\,\,$ the set of
its Fourier 
transforms. Note that in the classical sense
$$
\widehat{f}(\xi)=c_nH_{(n-2)/2}f_o(||\xi||), \quad f\in
L^p_{rad}({\bold R}^n),
\quad 1 \le p <2n/(n+1) . 
$$

By the convolution inequality there follows for $m\in[L^p_
{rad}(\bold R^n)]\,\,\widehat{}\,$ that
$$
T_m: L^1 \to L^p,  \qquad [T_m\varphi]\,\,\widehat{}\,\,(\xi): = 
m(||\xi||)\,\,\widehat{\varphi}(\xi), \qquad \varphi \in \Cal
S({\bold R}^n) 
$$
is bounded
and it is well known that if a bounded convolution operator from
$L^1({\bold 
R}^n)$ to $L^p({\bold R}^n)$, $1<p<\infty$, is generated by some
radial $m$ 
then $m\in [L^p_{rad}({\bold R}^n)]\,\,\widehat{}\,\,\,$. For
this instance we 
say that $m(||\xi||)\in M^{1,p}({\bold R}^n)$ and define 
$||m(||\xi||)||_{M^{1,p}({\bold R}^n)}$ to be the operator norm
of $T_m$
which is equal to the $L^p(\bold R^n, dx)$-norm of $H_{(n-2)/2}
m(||x||)$. 
Setting $\overline{\xi}=(\xi,\xi_{n+1})$ and
$\overline{\overline{\xi}}=(\overline{\xi},\xi_{n+2})$ we have
\proclaim{Corollary~3\.3}
Let  $1<p<q<\infty$ and $n\ge2$ be an integer. There holds

\noindent a)
$$
\bigg|\bigg|m(||\overline{\overline{\xi}}||)\bigg|\bigg|_{M^{1,q}
({\bold R}^{n+2})}
\le C\,\bigg|\bigg|\,||\xi||^2\,\,m(||\xi||)\bigg|\bigg|_{M^{1,p}
({\bold R}^n)}, \qquad  \frac1q 
=\frac1p - \frac{2}{(n+2)p}\,\,\,;
$$
b)
$$
\bigg|\bigg|m(||\xi||)\bigg|\bigg|_{M^{1,q}({\bold R}^n)}
\le
C\,\bigg|\bigg|m(||\overline{\xi}||)\bigg|\bigg|_{M^{1,p}({\bold
R}^{n+1})}\,,\qquad  \frac1q 
=\frac1p - \frac{1}{np'}\,\,\,.
$$
\endproclaim   
\demo{Proof}
For part a) choose
$\nu =(n-2)/2$ and $\mu =n/2$ for an integer $n\ge 2$ in
Theorem~3\.1.
The assumption there that $m$ is smooth  may be dropped since any
$H_\nu m(||x||)\in L^p(\bold R^n)$ can be
approximated in $L^p(\bold R^n)$ by smooth rapidly decreasing 
$H_\nu m_k(||x||)$ with $m_k\to m$ in $S'(\bold R^n)$, 
thus (3.11) gives the assertion a) for an arbitrary radial 
$m\in M^{1,p}({\bold R^n})$. Part b) follows similarly from
Theorem~3\.2 when choosing $\mu = (n-2)/2$ and $\nu =(n-1)/2$. 
\enddemo
{\bf Remarks}. 1) The results of Corollary~3.3 are best possible
for 
$1<p<q<2n/(n+1)$ in the following sense. For part a) consider the
example 
$$
m(t)=t^{-(n+2)/q'}(1+\log^2t)^{-1}=t^{-2}t^{-n/p'}(1+\log^2t)^{-1
}.
$$ 
By a criterion in [T] we have
$m(||\overline{\overline{\xi}}||)\in
[L^q({\bold R}^{n+2})]\,\,\widehat{}\,\,$ but $m(||\overline{
\overline{\xi}}||)$ does not belong to any
other space $[L^r({\bold R}^{n+2})]\,\,\widehat{}\,\,$, $r\neq
q$,
which follows directly by (3.12) on account of H\"older's
inequality.
The same reasoning applies to the right hand side.
Concerning Part b) we rewrite this example in the form
$$
m(t)=t^{-(n+1)/p'}(1+\log^2t)^{-1}=t^{-n/q'}(1+\log^2t)^{-1},
$$
and argue as in the case of Part a) which gives that also Part b) 

is best possible in the previous sense.

\noindent 2) 
Note that part a) is in the spirit of the following result due to
Coifman and Weiss, [CW, p\.33-45].
$$
\bigg|\bigg|m(||\overline{\overline{\xi}}||)\bigg|\bigg|_{M^{p,p}
(\bold R^{n+2})}\le
C\,\,\bigg|\bigg|n\,m(||\xi||)+||\xi||\,m'(||\xi||)\bigg|\bigg|
_{M^{p,p}(\bold R^n)}\,\,\,,
\tag3.13
$$
which is only good for $p$ near 1 or infinity. That the right
side of
(3\.13) contains an expression of type $tm'(t)$ is only natural
in view of the neccesary conditions for radial Fourier
multipliers in [GT, p\.412] (for $p$ near 1). 
These conditions also indicate that part a) of Corollary 3.3 is a
natural result ($p\le q<\frac {2n}{n+1})$; for the neccessary
conditions arising from the right side 
guarantee quite precisely the neccessary conditions arising from
the
left side.
Part b) of the corollary is in the spirit of  the well known
deLeeuw 
restriction result for Fourier multipliers (see e\.g\. [To], p\.
265) 
which by duality and the  Riesz interpolation theorem implies
$$
\bigg|\bigg|m(||\xi||)\bigg|\bigg|_{M^{q,q}({\bold R}^n)}
\le
C\,\bigg|\bigg|m(||\overline{\xi}||)\bigg|\bigg|_{M^{p,p}({\bold
R}^{n+1})}\,\,\,,
$$
$1\le\min\{p,p'\}\le q \le\max\{ p,p'\} \le \infty$.

\subhead
4. Density theorems
\endsubhead
In this section we prove density theorems which were announced
and
used in \S 2. Because they are of some independent interest and,
perhaps, could be used for other purposes, we prove these
theorems in
a more general form than we actually need them. The results we
obtain
are generalizations to the Hankel transform setting of density
theorems proved by Muckenhoupt, Wheeden and Young, [MWY].
In other words we extend the results of Section 2 of [MWY] from
the
cosine transform setting that corresponds to the case $\nu=-1/2$
to
general $\nu\ge -1/2$. Hence in what follows we restrict the
attention
to $\nu>-1/2$ only.
Needless to say we follow the ideas of [MWY] fairly closely.

If not otherwise stated the letter $k$ will always denote an
integer.
Recall that $\Cal S(\bold{R}_+)$ denotes the space of
restrictions to
$(0,\infty)$ of even Schwartz functions on $\bold{R}_+$ and
$C^\infty_o$ denotes the space of $C^\infty$ functions with
compact
support in $(0,\infty)$. Recall also that $H_\nu$ is a bijection
on
$\Cal S(\bold{R}_+)$. Observing that for even $f$ in $\Cal
S(\bold{R}_+)$ 
we have $f'(0)=0$ it is readily checked that the differential
operator
$L_\nu$ can be extended to even Schwartz functions by setting 
$L_\nu f(0)=2(\nu+1)f''(0)$. Thus,
if the powers of the operator $L_\nu$ are now defined in the
usual way:
$L_\nu^1=L_\nu$ and  $L_\nu^k=L_\nu(L_\nu^{k-1})$, $k>1$,
iterating the process we can regard $L_\nu^k f$, $k=0,1,\dots$,
to be a
function in $\Cal S(\bold{R}_+)$. 

\proclaim{Lemma~4.1} 
If $f\in\Cal S(\bold{R}_+)$ satisfies
$$
\int_0^\infty x^{2j}f(x)x^{2\nu+1}dx=0, \qquad j=0,1,\dots,k,
\tag4.1 
$$
then $(H_\nu f)^{(j)}(0)$, the derivatives of $H_\nu f$ at zero,
vanish for $j=0,1,\dots,2k$.
\endproclaim
\demo{Proof}
It follows from (1.8) that
$$
L^j_\nu H_\nu f(x)=(-1)^jH_{\nu}((\cdot)^{2j}f)(x), \qquad
x>0.\tag4.2
$$
Hence $L^j_\nu H_\nu f(0)=0$ for $j=0,1,\dots,k$ which implies 
$(H_\nu f)^{(j)}(0)=0$ for $j=0,2,\dots,2k$. It is obvious that
the
same holds for odd $j$'s. 
\enddemo
By $Q_k(\nu)$, $k\ge0$, we will denote the set of functions $f$
in
$L^2(dm_\nu)\cap L^{1,2k}(dm_\nu)$ that satisfy (4.1) and,
if $k<0$, we set $Q_{k}=L^2(dm_\nu)$; then we define
$C^\infty_o(k,\nu)
=C^\infty_o\cap Q_k(\nu)$.
\proclaim{Lemma~4.2} 
If $1\le p<\infty$, $\gamma>-1$,
$2k>-2+(\gamma+1)/p$ then every function $f$ in
$C^\infty_o(k,\nu)$
is approximated by functions from $H_{\nu}(C^\infty_0)$ in both
$L^{p,\gamma}(dx)$ and $L^2(dx)$ norms.
\endproclaim
\demo{Proof}
Let $\phi_n(x)$ be the sequence of functions on $(0,\infty)$
defined 
as in the proof of Lemma~6.2 in [MWY]:
$\phi_n(x)=\phi(nx)$ if $0<x\le 1/n$, $\phi_n(x)=\phi(x/n)$ if
$x\ge n$
and $\phi_n(x)=0$ if $n^{-1}\le x\le n$, where $\phi$ is a fixed
$C^\infty$
function on $(0,\infty)$ with $\phi(x)=0$ for $1/2\le x\le2$,
$\phi(x)=1$
for $0<x\le1/4$ and $x\ge4$ and $0\le\phi(x)\le1$ elsewhere.
Given $f$ in $C^\infty_o(k,\nu)$ define 
$f_n=H_\nu(H_\nu f\cdot(1-\phi_n))$. Since $1-\phi_n\in
C^\infty_o$
then $f_n\in H_\nu(C^\infty_o)$. The convergence of $f_n$ to $f$
in $L^2(dx)$ is immediate. To prove that $f_n$ approaches $f$ in 
$L^{p,\gamma}(dx)$ norm we write 
$$
||f-f_n||_{p,\gamma}\le||(f(x)-f_n(x))(1+x)^{2(k+1)}||_\infty
||(1+x)^{-2(k+1)}||_{p,\gamma}
$$
and note that the last norm on the right is finite due to
assumptions
on $p, \gamma$ and $k$. Moreover, by (4.2)
$$
\align
||(f-f_n)(1+x)^{2(k+1)}||_\infty
&\le C||f-f_n||_\infty+C||x^{2(k+1)}(f-f_n)||_\infty\\
&\le C||H_\nu(f-f_n)||_1+C||L_\nu^{k+1}(H_\nu f-H_\nu
f_n)||_1\,,
\endalign
$$
where $||\cdot||_1$ denotes the norm in $L^1(dm_\nu)$.
The fact that $H_\nu f-H_\nu f_n=H_\nu f\cdot\phi_n$ shows that 
$||H_\nu(f-f_n)||_1\to 0$ as $n\to\infty$. 
To estimate the remaining term we use the following Leibniz' rule
for
the $(k+1)$th power of the operator $L_\nu$
$$
L_\nu^{k+1}(H_\nu f\cdot\phi_n)=
\sum_{1\le i+j\le 2(k+1)}c_{ij}x^{-2(k+1)+i+j}(H_\nu
f)^{(i)}\phi_n^{(j)}.
$$
This may be proved by induction. We now consider the
$L^1(dm_\nu)$ norm of each summand in the sum above separately.
Fixing $i,j$, $1\le i+j\le2(k+1)$ we have to show that the
quantities
$$
n^{-j}\int_n^\infty|(H_\nu f)^{(i)}(x)\phi^{(j)}(\frac
xn)|x^{-2(k+1)+i+j+ 2\nu+1}dx \tag4.3
$$ 
and
$$
n^{j}\int_0^{1/n}|(H_\nu f)^{(i)}(x)\phi^{(j)}(xn)|x^
{-2(k+1)+i+j+ 2\nu+1}dx \tag4.4
$$ 
tend to 0 as $n\to\infty$. This is easily seen for (4.3) since
$\phi^{(j)}$ is bounded and $(H_\nu f)^{(i)}$ is of rapid
decrease at
$\infty$. For (4.4), consider first the case $i=2(k+1)$. Then
$j=0$
and (4.4) is bounded by $Cn^{-(2\nu+2)}$. If $0\le i\le 2k+1$
then 
by Taylor's formula and Lemma~4.1 the estimate $|(H_\nu
f)^{(i)}(x)| \le Cx^{2k+1-i}$ follows. This shows that (4.4) is
bounded by
$Cn^{-(2\nu+1)}$ and finishes the proof of Lemma~4.2.
\enddemo

\proclaim{Lemma~4.3} 
If $1\le p<\infty$, $\gamma>-1$, then every function $f$ in 
$Q_k(\nu)\cap L^{p,\gamma}(dx)$
is approximated by functions from $C^\infty_0(k,\nu)$ in both
$L^{p,\gamma}(dx)$ and $L^2(dx)$ norms.
\endproclaim
The proof of Lemma~4.3, with minor changes, is the same as the
proof 
of Lemma~6.6 in [MWY]. Let us mention at this point that for our
future purposes we will use, for given $k$, a sequence of
$C^\infty$
functions $\{\alpha_j(x)\}_0^k$, the same as in Lemma~6.5 of
[MWY]
except for the fact that their supports are separated from zero,
say,
they are contained in $1/4\le x\le3/4$. It can be checked that
this
requirement is not essential. Recall that an important feature of
$\alpha_j$'s is the fact that
$$
\int_0^\infty x^i\alpha_j(x)dx=\delta_{i,j}, \tag4.5
$$ 
$0\le i,j\le k$, where $\delta_{i,j}$ is the Kronecker delta.

\proclaim{Lemma~4.4} 
If $1\le p<\infty$, $\gamma>-1$, $2k<-1+(\gamma+1)/p-(2\nu+1)$
then every function $f$ in $C^\infty_o$
is approximated by functions from $C^\infty_0(k,\nu)$ in 
$L^{p,\gamma}(dx)$ norm.
\endproclaim
\demo{Proof}
If $k$ is negative the statement is obvious. Let $k\ge0$ and
take $\alpha_0, \alpha_1,\dots,\alpha_{2k}$ satisfying (4.5) for
$0\le
i,j\le2k$, supported in $1/4\le x\le3/4$ and, given $f\in
C^\infty_o$,
define
$$
f_n(x)=f(x)-\sum_{i=0}^kn^{2i+1}\alpha_{2i}(nx)x^{-(2\nu+1)}
\int_0^\infty f(t)t^{2i}t^{2\nu+1}dt\,.
$$ 
Then $f_n\in C^\infty_o(k,\nu)$ and the required convergence
$f_n\to
f$, $n\to\infty$, holds in $L^{p,\gamma}(dx)$.
\enddemo
\proclaim{Lemma~4.5} 
If $1\le p<\infty$, $\gamma>-1$, $2k>-3+(\gamma+1)/p-(2\nu+1)$
then every function $f$ in $Q_k(\nu)\cap L^{p,\gamma}(dx) $
is approximated by functions from $C^\infty_0(k+1,\nu)$ in both
$L^{p,\gamma}(dx)$ and $L^2(dx)$ norms.
\endproclaim
\demo{Proof} By Lemma~4.3 we can assume that $f$ is in
$C^\infty_o(k,\nu)$.
Again the statement is obvious if $k<-2$. Hence, assume $k\ge-1$
and take $\alpha_0, \alpha_1,\dots,\alpha_{2(k+1)}$ satisfying
(4.5) for $0\le
i,j\le2(k+1)$ and define
$$
f_n(x)=f(x)-n^{-2(k+1)-1}\alpha_{2(k+1)}(x/n)x^{-(2\nu+1)}
\int_0^\infty f(t)t^{2(k+1)}t^{2\nu+1}dt\,.
$$ 
Then $f_n$ is in $C^\infty_0(k+1,\nu)$ and $f_n$ converges to $f$
in
$L^{p,\gamma}(dx)$ and $L^2(dx)$.
\enddemo
\proclaim{Lemma~4.6} 
If $1<p<\infty$, $\gamma>-1$, $2k=-3+(\gamma+1)/p-(2\nu+1)$
then every function $f$ in $Q_k(\nu)\cap L^{p,\gamma}(dx) $
is approximated by functions from $Q_{k+1}(\nu)\cap
L^{p,\gamma}(dx)$ in both
$L^{p,\gamma}(dx)$ and $L^2(dx)$ norms.
\endproclaim
\demo{Proof} We can consider the case $k\ge-1$ only. Take
$\alpha_0, 
\alpha_1,\dots,\alpha_{2(k+1)}$ satisfying (4.5) for $0\le
i,j\le2(k+1)$ and define
$$
g_n(x)=\frac{\chi_{[e,n]}(x)x^{-(2\nu+1)}}{x^{2(k+1)+1}\log
x\log\log
n}
-\sum_{i=0}^{k}\frac{\alpha_{2i}(x)x^{-(2\nu+1)}}{\log\log n}
\int_e^n \frac{t^{2i-2(k+1)-1}}{\log t}dt\,.
$$ 
Then $\int_0^\infty g_n(x)x^{2i}x^{2\nu+1}dx$ equals 0 for
$i=0,1,\dots,k$ and is 1 for $i=k+1$. An  argument shows that
$g_n\to 0$ in $L^2(dx)$ and convergence of $g_n$ to 0 in
$L^{p,\gamma}(dx)$ is implied by the fact that $\int_e^\infty
(x(\log x)^p)^{-1}dx<\infty$ for $p>1$.
By Lemma~4.3 we can assume $f$ to be in $C^\infty_0(k+1,\nu)$.
Define
$$
f_n(x)=f(x)-g_n(x)\int_0^\infty f(t)t^{2(k+1)}t^{2\nu+1}dt\,.
$$
Then $f_n\in C^\infty_0(k+1,\nu)$ and the properties of $g_n$
imply
the desired convergence for $f_n$.
\enddemo
\proclaim{Theorem~4.7} 
If $1<p<\infty$, $\gamma>-1$
then $H_\nu(C_o^\infty)$ is dense in $L^{p,\gamma}(dx) $.
If, in addition $\gamma$ is not of the form $\gamma=2k+2\nu+1$
then $H_\nu(C_o^\infty)$ is dense in $L^{1,\gamma}(dx) $.
\endproclaim
\demo{Proof} 
Fix $p$ and $\gamma$ and choose $k$ to be an integer satisfying
$-3+(\gamma+1)/p-(2\nu+1)\le 2k<-1+(\gamma+1)/p-(2\nu+1) $
if $p>1$ and $-2+\gamma-(2\nu+1)<2k<\gamma-(2\nu+1)$ if  $p=1$.
Since $C_o^\infty$ is dense in $L^{p,\gamma}(dx) $ it is
sufficient
to approximate functions from $C_o^\infty$ only. Lemma~4.4 allows
further to restrict the attention to functions from
$C_o^\infty(k,\nu)$. By applying Lemma~4.5 or Lemma~4.6 several
times
and then applying Lemma~4.3, if necessary, we conclude that every
function from $C_o^\infty(k,\nu)$ is approximated by functions
from 
$C_o^\infty(k+r,\nu)$, where $r$ is a positive integer such that
$2(k+r)>-2+(\gamma+1)/p$. Using Lemma~4.2 finishes the proof of
the
theorem. 
\enddemo
\proclaim{Corollary~4.8} 
If $1<p<\infty$, $\beta>-1-p(\nu+1/2)$ then $\Cal
H_\nu(C_o^\infty)$ is dense in $L^{p,\beta}(dx)$.
If, in addition $\beta$ is not of the form $\beta=2k+\nu+1/2$
then $\Cal H_\nu(C_o^\infty)$ is dense in $L^{1,\beta}(dx) $.
\endproclaim
\demo{Proof} 
The corollary follows from Theorem~4.7 by using (1.4), the fact
that
one has $x^{\nu+1/2}C_o^\infty=C_o^\infty$ and the remark that
multiplication
by $x^{\nu+1/2}$ is an isometric bijection between
$L^{p,\gamma}(dx)$
and $L^{p,\gamma-p(\nu+1/2)}(dx)$.
\enddemo

\enddocument